\input amstex
\documentstyle {amsppt}

\pageheight{50.5pc} \pagewidth{32pc}

\topmatter
\title{}
A uniqueness theorem for the martingale problem describing a
diffusion in  media with membranes
\endtitle
\rightheadtext{A uniqueness theorem for the martingale problem}
\leftheadtext{Olga V. Aryasova, Mykola I. Portenko}

\author
Olga V. Aryasova, Mykola I. Portenko
\endauthor
\address
Institute of Geophysics, National Academy of Sciences of Ukraine,
Palladina pr. 32, 03680, Kiev-142, Ukraine
\endaddress
\email oaryasova\@mail.ru
\endemail
\keywords Diffusion process, martingale problem, uniqueness of
solution
\endkeywords

\address
Institute of Mathematics,  National Academy of Sciences of
Ukraine, Tereshchenkivska str. 3, 01601, Kiev, Ukraine
\endaddress
\email  portenko\@imath.kiev.ua
\endemail
\subjclass 60J60, 60J35
\endsubjclass
\abstract We formulate a martingale problem that describes a
diffusion process in a multidimensional Euclidean space with a
membrane located on  a given smooth surface and having the
properties of skewing and delaying. The theorem on the existence
of no more than one solution to the problem is proved.
\endabstract
\endtopmatter

\document
\head { Introduction}
\endhead

Let $S$ be a  given closed bounded surface in ${\Bbb R}^d$ that
divides the space ${\Bbb R}^d$ into two open parts: the interior
domain $D^i$ and the exterior one $D^e$, $D$ is the union of them.
The surface $S$ is assumed to be smooth enough (see Section 1 for
the precise assumptions) so that at any point of $S$ there is a
well-defined normal. We denote by $\nu(x)$ for $x\in S$ the unit
vector of outward normal to $S$ at the point $x$. Let $A(x),\ x\in
S,$ be a given real-valued continuous function and, for each $y\in
{\Bbb R}^d$, let $b(y)$ be a symmetric positively definite linear
operator in $\Bbb R^d$. The function $(b(y))_{y\in \Bbb R^d}$ is
supposed to be bounded and H\"older continuous. For $x\in S,$ the
vector $N(x)=b(x)\nu(x)$ is called the co-normal vector to $S$ at
the point $x.$  Consider the stochastic differential equation in
${\Bbb R}^d$
$$
dx(t)=A(x(t))N(x(t))\hbox{1\!\!\!\!\;\,{I}}_S(x(t))dt+b(x(t))^{1/2}\hbox{1\!\!\!\!\;\,{I}}_D(x(t))dw(t),
\tag1
$$
where $(w(t))_{t\geq0}$  is a standard Wiener process in ${\Bbb
R}^d,\ \ \hbox{1\!\!\!\!\;\,{I}}_\Gamma$ is  the indicator
function of a set $\Gamma \subset {\Bbb R}^d$. As was shown in
[1],  this equation has infinitely many solutions. Consequently,
if a solution to (1) is treated as that to the corresponding
martingale problem, the latter turns out not to be well-posed.

Each solution constructed in [1] is determined by a representation
of the function $A(x),\ x \in S,$ in the form
$A(x)=\frac{q(x)}{r(x)}$, where $q(\cdot)$ and $r(\cdot)$ are
continuous functions on $S$ taking their values in $[-1,1]$ and
$(0,+\infty)$, respectively. Thus, the formulation of the
well-posed martingale problem must involve these functions.

A solution to (1) was constructed in [1] as a continuous Markov
process $(x(t))_{t\geq0}$ in ${\Bbb R}^d$ obtained from a
$d$-dimensional diffusion process with its diffusion operator
$b(\cdot)$ and zero drift vector by two transformations. The first
transformation is skewing the diffusion process on $S$. The skew
is determined by the function $q(\cdot)$. As a result, one get a
continuous Markov process $(x_0(t))_{t\geq0}$ in ${\Bbb R}^d$ such
that its trajectories satisfy the following stochastic
differential equation (see [2], Ch. 3)
$$
dx_0(t)=q(x_0(t))\delta_S(x_0(t))N(x_0(t))dt+b(x_0(t))^{1/2}dw(t),
\tag2
$$
where $(\delta_S(x))_{x\in{\Bbb R}^d}$   is a generalized function
on ${\Bbb R}^d$  that acts on a test function
$({\varphi}(x))_{x\in{\Bbb R}^d}$ according to the following rule
$$
{\langle}\delta_S, {\varphi}{\rangle}=\int_S{\varphi}(x)d\sigma
$$
(the integral in this equality is a surface integral).

To do the second transformation determined by a given function
$r(\cdot):S\to(0,\infty)$, one should put for $t\geq0$
$$
\zeta_t=\inf\bigg\{s: \
s+\int^s_0r(x_0(\tau))\delta_S(x_0(\tau))d\tau\geq t \bigg\}
$$
and define
$$ x(t)=x_0(\zeta_t),\ t\geq 0.
$$
Here the functional
$$
\eta_t=\int^t_0r(x_0(\tau))\delta_S(x_0(\tau))d\tau, \ t\geq0,
$$
of the process $(x_0(t))_{t\geq0}$ is well defined as an additive
homogeneous continuous functional (see [2], Ch.3).  As is known
(see [3], Theorem 10.11), the process $(x(t))_{t\geq0}$   is a
continuous Markov process in ${\Bbb R}^d$  as a result  of the
random change of time for the process $(x_0(t))_{t\geq0}.$

The following  observation gives us a suggestion how to formulate
correctly the martingale problem for the process $(x(t))_{t\geq0}$
corresponding to a given pair of functions $q(\cdot)$ and
$r(\cdot)$. Namely, fix  an orthonormal basis  in ${\Bbb R}^d$ and
denote by $x_j$    for $j=1,2,\ldots, d$   the coordinates of a
vector $x\in{\Bbb R}^d$ and by $b_{jk}(x)$   for $j,k=1,2,\ldots,
d$  the elements of the matrix of the operator $b(x)$ in that
basis. For  a given continuous bounded function $\varphi$ on $\Bbb
R^d$ with real values, we put $u(t,x,\varphi)=\Bbb
E_x\varphi(x(t)),\ t\geq0$ and $x\in R^d$. Then this function is
continuous in the arguments $t\geq0$ and $x\in \Bbb R^d$ and turns
out to satisfy the following conditions:

1) it satisfies the equation
$$
\frac{{\partial} u}{{\partial}
t}=\frac{1}{2}\sum^d_{i,j=1}b_{ij}(x)
\frac{{\partial}^2u}{{\partial} x_i{\partial} x_j}\tag3
$$
in the domain $t>0, \ x \in D$;

2) it satisfies the equation
$$ r(x)\frac{\partial u}{\partial
t}=\frac{1+q(x)}{2}\frac{\partial u(t,x+)}{\partial N(x)} -
\frac{1-q(x)}{2}\frac{\partial u(t,x-)}{\partial N(x)}
$$
for  $t>0,\ x\in S$;

3)  the initial condition
$$
u(0+,x)=\varphi(x)
$$
is held for all $x\in \Bbb R^d$.

In Section 1 we give a correct form of the martingale problem
desired. Our aim is to show that the solution to that problem is
unique. We obtain the statement from the uniqueness theorem for
the boundary process (see Section 2 for the precise definition) by
the Strook-Varadhan method from [4]. The particular case of an
identity diffusion matrix and $S$ being a hyperplane was
investigated in [5]. One can also find there some further
discussion of the topic.

\head {1. The martingale problem}
\endhead
From now on we assume that, for each $x\in{\Bbb R}^d$,
$b(x)=\left(b_{ij}(x)\right)_{i,j=1}^d$ is a symmetric $d\times
d$-matrix satisfying the following conditions  which we call the
conditions $J$

1) there are two positive constants $C_1$ and $C_2,\ 0<C_1\leq
C_2,$ such that
$$
C_1|\theta|^2\leq(b(x)\theta,\theta)\leq C_2|\theta|^2
$$
is valid for all $\theta,x\in{\Bbb R}^d$.

2) for all $x,x'\in{\Bbb R}^d, \ i,j=1,2,\dots,d$,
$$
|b_{ij}(x)-b_{ij}(x')|\leq L|x-x'|^\alpha,\tag4
$$
where $L$ and $\alpha$ are positive constants, $\alpha\leq 1$.

Suppose $S$ belongs to the class $H^{2+\varkappa}$ for some
$\varkappa\in (0,1)$ (see [6], Ch. 4, \S \ 4). By $\delta$ we
denote the minimal one of the numbers $\alpha$ from (4) and
$\varkappa$.

Suppose a continuous function $q(\cdot): S\to[-1,1]$ and a
continuous bounded function $r(\cdot): S\to[0,+\infty)$ are fixed.

$\Omega$ stands for the space of all continuous ${\Bbb
R}^d$-valued  functions on $[0,+\infty)$, $\Cal M_t$ denotes the
$\sigma$-algebra generated by  $x(u)$  for $0\leq u \leq t$. If
$t=\infty$,  $\Cal M_t$ will be denoted by $\Cal M$.

We say that a function $f$ belongs to the class $F$ if

1)$\ f$ is continuous and bounded in $(t,x)$ on
$[0,+\infty)\times{\Bbb R}^d$;

2) $f$ has a  continuous and bounded derivative in $t$ on
$[0,+\infty)\times{\Bbb R}^d$;

3) $f$ has continuous and bounded derivatives in $x$ on
$[0,+\infty)\times D$ up to the second order;

4) for all $t\in [0,+\infty)$ and $x\in S$ there exist the
non-tangent limits $ \dfrac{\partial f(t,x+)}{\partial N(x)}$ and
$\dfrac{\partial f(t,x-)}{\partial N(x)}$ from the side $D^e$ and
$D^i$, respectively, and the function
$$
Kf(t,x)= \frac{1+q(x)}{2}\frac{\partial f(t,x+)}{\partial N(x)} -
\frac{1-q(x)}{2}\frac{\partial f(t,x-)}{\partial N(x)}
$$
is continuous and bounded on $[0,+\infty)\times S$.
\definition
{ Definition 1} Given $x\in{\Bbb R}^d,$ a probability measure
${\Bbb P}_x$ on $\Cal M$ is a solution to the submartingale
problem starting from $x$ if

1) ${\Bbb P}_x\{x(0)=x\}=1;$

2) the process
$$
X_f(t)=f(t,x(t))-\int^t_0\hbox{1\!\!\!\!\;\,{I}}_D(x(u))\left(\frac{\partial
f}{\partial
u}+\frac{1}{2}\sum_{i,j=1}^{d}b_{ij}(x(u))\dfrac{\partial^2
f}{\partial x_i \partial x_j}\right)(u,x(u))du,\ t\geq0,
$$
is a ${\Bbb P}_x$-submartingale whenever $f$ belongs to $F$ and
satisfies the inequality
$$
r(x)\frac {\partial f(t,x)}{\partial t}+Kf(t,x)\geq0 \ \hbox{for}
\ t\geq0 \hbox{ and} \ x\in S.
$$
\enddefinition
\remark {Remark 1} One can verify that the transition probability
of the process $(x(t))_{t\geq0}$ described in the Introduction is
a solution to the submartingale problem ([1]).
\endremark
\vskip 5 pt
Define the function $\phi$ on ${\Bbb R}^d$ by the
equality $\phi(x)=d(x,S):=\inf\{d(x,y): y\in S\}$, where
$d(\cdot,\cdot)$ is the Euclidean metric on ${\Bbb R}^d$ . Then

1) $S=\{x\in{\Bbb R}^d: \phi(x)=0\},$ $D=\{x\in{\Bbb R}^d:
\phi(x)>0\},$

2) $K\phi(x)\equiv1$ on $S.$ \vskip5pt \remark {Remark 2} The
function $\phi$  does not belong to the class $F$ because of its
unboundedness. To overcome this, we choose, for each $m\geq1$, a
non-increasing infinitely differentiable function $\eta_m$ defined
on $[0,\infty)$ and having compact support such that
$0\leq\eta_m\leq1,\ \eta_m\equiv 1$ on $[0,m]$, $\eta_m\equiv 0$
off $[0,m+1]$ and the derivatives of $\eta_m$ up to the second
order are uniformly bounded. Set $\phi_m(x)=\eta_m(d(x,S)) \!
\cdot \! \phi(x),\ x\in {\Bbb R}^d$. Then $\phi_m$ belongs to $F$.
Hence $X_{\phi_m}$ is a ${\Bbb P}_x$-submartingale. Clearly,
$\phi_m(x)\rightarrow \phi(x)$ monotonically as $m\rightarrow
\infty$ and $\sum_{i,j=1}^{d}b_{ij}(x)\dfrac{\partial^2
\phi_m(x)}{\partial x_i
\partial x_j}$ tends to 0 boundedly. So $X_{\phi}(t)=\phi(x(t))$
is a ${\Bbb P}_x$-submartingale.
\endremark

The following proposition gives a reformulation of the
submartingale problem into a martingale one. \proclaim
{Proposition 1} Given $x\in{\Bbb R}^d,$ the probability measure
${\Bbb P}_x$ on $\Cal M$ solves the submartingale problem starting
from $x$ iff \ ${\Bbb P}_x\{x(0)=x\}=1$ and there exists a
continuous non-decreasing $(\Cal M_t)$-adapted process $\gamma(t),
\ t\geq0,$ such that

1) $\gamma(0)=0,$ $\Bbb E\gamma(t)<+\infty$  for all $t\geq0;$

2) {\rm
$\gamma(t)=\int^t_0\hbox{\hbox{1\!\!\!\!\;\,{I}}}_S(x(u))d\gamma(u),
t\geq0;$}

3) the process {\rm $$
f(t,x(t))-\int^t_0\hbox{1\!\!\!\!\;\,{I}}_D(x(u))
\left(\frac{\partial f}{\partial
u}+\frac{1}{2}\sum_{i,j=1}^{d}b_{ij}(x(u))\dfrac{\partial^2
f}{\partial x_i \partial x_j}\right)(u,x(u))du-
$$
}
$$
-\int^t_0(r\frac{\partial f}{\partial u}+Kf)(u,x(u))d\gamma(u),
t\geq0,
$$
is a ${\Bbb P}_x$-martingale for any $f$ belonging to $F.$

If \ ${\Bbb P}_x$  is such a solution, then $\gamma(t)$ is
uniquely determined, up to ${\Bbb P}_x$-equivalence, by the
condition that
$$
\phi(x(t))-\gamma(t),\ t\geq0,
$$
is a ${\Bbb P}_x$-martingale.
\endproclaim
\demo { Proof} The existence of a solution to this problem was
established  in [7]. The proof of the last statement is similar to
that of Theorem 2.5 in [4].
\enddemo
\proclaim { Corollary 1} For each $x\in{\Bbb R}^d, \ t\geq0,$ the
equality {\rm
$$
\int^t_0\hbox{1\!\!\!\!\;\,{I}}_S(x(u))du=\int^t_0r(x(u))d\gamma(u)
$$
} is held ${\Bbb P}_x$-almost surely.
\endproclaim
\proclaim
 { Corollary 2}  If $x\in S,$ then
${\Bbb P}_x\{\gamma(t)>0, t>0\}=1.$
\endproclaim

These assertions can be verified like Corollaries 1,2 in [5].

\head {2. A uniqueness theorem for a boundary process}
\endhead

Let ${\Bbb P}_x$  be a solution to the submartingale problem
starting from $x\in S.$ Then there exists a process $\gamma(t), \
t\geq0,$ that has the properties stated in Proposition 1. For
$\theta\geq0$, we put $\tau(\theta)=\sup\{t\geq0:
\gamma(t)\leq\theta\}.$ Define
$T(\omega)=\lim_{t\to+\infty}\gamma(t)$. Assume, that
$T(\omega)=+\infty$ a.s. Then the process
$y(\theta)=x(\tau(\theta))$ is defined for all
$0\leq\theta<\infty$. It is not hard to see that the process
$\tau(\theta)$  and, consequently,  the process $y(\theta)$ are
right continuous processes having no discontinuities of the second
kind, and the latter takes on its values on $S$. Since the
starting point is on $S$  we have $\gamma(t)>0$   for $t>0$ almost
surely, i.e. $\tau(0)=0$ and $y(0)=x.$ Following Strook and
Varadhan [4] we define the $(d+1)$-dimensional process
$(\tau(\theta), y(\theta)),\  \theta\geq0,$  and call it the
boundary process starting from $x.$ If $T(\omega)<\infty$ with
positive probability we put $(\tau(\theta), y(\theta))=\infty$ for
$\theta\geq T(\omega)$.

Further on we denote, by $C^{1,2}_0([0,+\infty)\times S)$,  the
class of functions on $[0,+\infty)\times S$ that have compact
supports with respect to $t$ and together with their first
$t$-derivative and two $x$-derivatives are continuous and bounded,
$C^\infty_0([0,+\infty)\times S)$  stands for the class of
infinitely differentiable functions on $[0,+\infty)\times S$
having compact supports with respect to $t$.

\proclaim { Proposition 2} For each $h\in
C^{1,2}_0([0,+\infty)\times S)$, there exists a function $Hh$ such
that

(i) it  belongs to the class $F$;

(ii) it is a solution to the equation
$$
\frac{\partial U}{\partial t}+\frac{1}{2}\sum^d_{i,j=1}b_{ij}(x)
\frac{{\partial}^2U}{{\partial} x^i{\partial} x^j}=0  \tag5
$$
on both   $[0,+\infty)\times D^i$ and $[0,+\infty)\times D^e$;

(iii) the relations
$$
Hh(t,x+)=h(t,x),     \tag 6
$$
$$
Hh(t,x-)=h(t,x),        \tag 7
$$
hold true for all $t\geq0$
and $x\in S.$
\endproclaim

\demo{ Proof} Assume that $\Theta$ is a domaine in $\Bbb R^d$. Set
$\Theta_T=(0,T)\times\Theta$ and denote, by $\overline{\Theta}_T$,
its closure. Let $H^{\delta/2+1,\ \delta+2}(\overline{\Theta}_T)$
be a corresponding H\"older space (see [6]), $H_T^{\delta/2+1,\
\delta+2}(\overline{\Theta}_T)$ stands for the set of all
functions from $H^{\delta/2+1,\ \delta+2}(\overline{\Theta}_T)$
which together with their first derivatives in $t$ are equal to
zero at the point $t=T$.
Notice that for all $T>0$, $h\in H^{\delta/2+1,\
\delta+2}(\overline{S}_T)$. Besides, there exists $T_0>0$ such
that $h=0$ if $t\geq T_0$. Therefore, $h\in H^{\delta/2+1,\
\delta+2}_{T_0}(\overline{S}_{T_0})$. By analogy to Theorem 5.2 in
[6] there is a uniquely defined function $h_{T_0}^i\in
H_{T_0}^{\delta/2+1,\ \delta+2}(\overline{D}^{\ i}_{T_0})$ which
satisfy equation (5) in the domains $(0,T_0)\times D^i_{T_0}$ and
boundary conditions (6). Remark that if $h(t,x)\equiv
 0$ on  $[A,B]\times S,$ where $A$ and $B$ are some constant, $0<A<B$, then
 the function equal to zero on $[A,B]\times \overline {D}^{\ i}$ is
 the unique solution to  problem (5), (6) belonging to $H_{B}^{\delta/2+1,\ \delta+2}
 ([A,B]\times
 \overline{D}^{\ i}_{B}).$ Similarly, there is a uniquely defined
 function  $h_{T_0}^e\in H_{T_0}^{\delta/2+1,\
 \delta+2}(\overline{D}^{\ e}_{T_0})$ that is a solution to (5), (7) in  $(0,T_0)\times
 D^e_{T_0}$. Now we define the function
 $$
 Hh=\cases
 h_{T_0}^i \ \hbox{on} \ \overline{D}^{\ i}_{T_0},\\
h_{T_0}^e \ \hbox{on} \ \overline{D}^{\ e}_{T_0}, \\
0, \ \hbox {otherwise}.
\endcases
$$
The function $Hh$ has all the required properties and this
completes the proof.
\enddemo

\vskip 10 pt

 Proposition 2 implies that for each $h\in C^{1,2}_0([0,+\infty)\times S)$,
$$
(\widetilde{K}h)(t,x)=r(x)\frac{\partial(Hh)(t,x)}{\partial t}+
$$
$$
+ \left[ \frac{1+q(x)}{2}\frac{\partial(Hh)(t,x+)}{\partial N(x)}-
\frac{1-q(x)}{2}\frac{\partial(Hf)(t,x-)}{\partial N(x)}\right]
$$
is well defined as a continuous and bounded function on
$[0,+\infty)\times S.$

\proclaim { Proposition 3} Suppose a probability measure ${\Bbb
P}_x$ solves the submartingale problem starting from $x\in S.$
Then the relation ${\Bbb P}_x\{(\tau(0), y(0))=(0,x)\}=1$ is held,
and for any function  $h\in C^{1,2}_0([0,+\infty)\times S)$ the
process
$$
h(\tau(\theta), y(\theta))-\int^\theta_0(\widetilde{K}h)(\tau(u),
y(u))du, \ \theta\geq 0,
$$
is a ${\Bbb P}_x$-martingale with respect to the filtration
$\left(\Cal M_{\tau(\theta)}\right)_{\theta\geq 0}$.
\endproclaim
\demo
 { Proof} The proof follows that of
Theorem 4.1 in [4].
\enddemo

Denote, by $\Cal D([0,+\infty),[0,+\infty)\times S)$, the class of
$[0,+\infty)\times S$-valued right-continuous functions on
$[0,+\infty)$ with no discontinuities of the second kind.
\definition
 { Definition 2} The uniqueness theorem is
valid for the boundary process if, for any given $x\in S$, there
is only one probability measure $\Bbb Q_x$,      on the space
$\Cal D([0,+\infty), [0,+\infty)\times S)$ such that

1) $\Bbb Q_x\{\tau(0)=0, y(0)=x\}=1;$

2) $(\tau(\theta), y(\theta))=\infty$ if $\theta>T(\omega)$;

3) for any function $h \in C^{1,2}_0([0,+\infty)\times S)$, the
process
$$
h(\tau(\theta), y(\theta)
)-\int^\theta_0(\widetilde{K}h)(\tau(u),y(u))du,\ \theta\geq0,
$$
is a $\Bbb Q_x$-martingale relative to  the natural
$\sigma$-algebras $\widetilde{\Cal M}_{\theta},\ {\theta\geq 0},$
in \break $\Cal D([0,+\infty),[0,+\infty)\times S).$
\enddefinition

To prove the uniqueness theorem for the boundary process, we will
make use of the following lemma.

\proclaim {Lemma} For each $\lambda>0, \ \psi\in
C^\infty_0([0,+\infty)\times S)$, the equation
$$
\lambda f-\widetilde{K}f=\psi\tag8
$$
has a unique solution in the class of all continuous functions  on
$[0,\infty)\times S$ having  compact supports with respect to $t$
and such that there exists a function $Hf$ on $[0,\infty)\times
\Bbb R^d$ satisfying conditions (i)-(iii) of Proposition 2.
\endproclaim
\demo {Proof}  We first prove the Lemma in the case of $r$ being
identically equal to $0$ on $S$.

Let $g_0(t,x,y), \ t>0,\ x\in \Bbb R^d, \ y\in\Bbb R^d$, be the
fundamental solution to the equation (3)(see [8], Ch. I). The
process $(x_0(t))_{t\geq0}$ solving equation (2) possesses a
transition probability density. We denote it by $G_0(t,x,y),\
t>0,\ x\in \Bbb R^d,\ y\in \Bbb R^d.$ As is proved in [1], for
$t>0,\ x\in{\Bbb R}^d,$   and $y\in{\Bbb R}^d,$ the representation
$$
G_0(t, x, y)=g_0(t, x, y)+\int^t_0d\tau\int_S\widetilde{V}(\tau,
x, z) \frac{{\partial} g_0(t-\tau, z, y)}{{\partial} N(z)}
q(z)d\sigma_z
$$ takes place,   where for $t>0,\ x\in {\Bbb R}^d,$    and $y\in S,\ \widetilde{V}(t,x,y)$  is the
solution to the following integral equation
$$
\widetilde{V}(t, x, y)=g_0(t, x,
y)+\int^t_0d\tau\int_S\widetilde{V}(\tau, x, z) \frac{{\partial}
g_0(t-\tau, z, y)}{{\partial} N(z)} q(z)d\sigma_z.
$$

Besides, the equality
$$
\widetilde{V}(t, x,
y)=\frac{1}{2}\left[G_0(t,x,y+)+G_0(t,x,y-)\right]\tag9
$$
is valid for $t>0,\ x\in {\Bbb R}^d,$    and $y\in S,$ and the
equality
$$
\frac{1+q(x)}{2}\frac{\partial G_0(t,x+,y)}{\partial N(x)} -
\frac{1-q(x)}{2}\frac{\partial G_0(t,x-,y)}{\partial N(x)}=0\tag10
$$
is held for $t>0,\ x\in S,$    and $y\in \Bbb R^d.$

For $\lambda\geq0,$ we define a function $G_\lambda$  of the
arguments $t>0, x\in{\Bbb R}^d,$  and $y\in{\Bbb R}^d$ by the
relation
$$
E_x({\varphi}(x_0(t))\exp\{-\lambda\eta_t\}) =\int_{{\Bbb
R}^d}{\varphi}(y)G_\lambda (t, x, y)dy
$$
that must be fulfilled for all $t>0,\ x\in{\Bbb R}^d,$ and
${\varphi}$ being a bounded measurable function on $\Bbb R^d.$
Then (see [1]) such a function exists and it can be found as a
solution to the pair of equations
$$
G_\lambda(t,x,y)=G_0(t,x,y)-\lambda\int^t_0d\tau\int_S\widetilde{V}(\tau,x,z)
G_\lambda(t-\tau, z,y)r(z)d\sigma_z,\tag11
$$
$$
G_\lambda(t,x,y)=G_0(t,x,y)-\lambda\int^t_0d\tau\int_S
G_\lambda(\tau,x,z) G_0(t-\tau, z,y)r(z)d\sigma_z\tag12
$$
in the domain $t>0,\ x\in {\Bbb R}^d,$    and $y\in \Bbb R^d$.
Moreover, there is no more than one solution to these equations
satisfying the inequality
$$
G_\lambda(t,x,y)\leq G_0(t,x,y).\tag13
$$

Then for all $t\geq 0 ,\ x\in \Bbb R^d,$ we can define the
function
$$
V_\lambda(t,x)=\int^\infty_td\tau\int_SG_\lambda^1(\tau-t,x,y)\psi(\tau,y)
d\sigma_y,
$$
where $G^1_\lambda(t,x,y)$ is the solution to the pair of
equations (11), (12) satisfying inequality (13) for $r(x)\equiv1$.

Notice, that $G_\lambda$ as a function of the third argument has a
jump at the points of $S$.
 Namely, for $t>0, x\in{\Bbb R}^d,$  and
$y\in S,$ the equations (11)  can be rewritten as follows
$$
G_\lambda(t,x,y)=\widetilde{V}(t,x,y)-\lambda\int^t_0d\tau\int_S
\widetilde{V}(\tau,x,z) G_\lambda(t-\tau, z, y)r(z)d\sigma_z.
$$

From this we can write the following relation for the function
$V_\lambda(t,x):$
$$
V_\lambda(t,x)= \int^\infty_td\tau\int_S \widetilde
V(\tau-t,x,y)\psi(\tau,y) d\sigma_y-\lambda \int^\infty_t
d\tau\int_S \widetilde {V}(\tau-t,x,y)V_\lambda(\tau,y) d\sigma_y.
\tag14
$$
The function $V_\lambda(t,x)$ has the following properties

1) $V_\lambda(t,x)$ satisfies conditions (i), (ii) of Proposition
2;

2) the equality
$$
\lambda V_\lambda(t,x)- \left[ \frac{1+q(x)}{2}\frac{\partial
V_\lambda(t,x+)}{\partial N(x)}- \frac{1-q(x)}{2}\frac{\partial
V_\lambda(t,x-)}{\partial N(x)} \right]=\psi(t,x)\tag15
$$
is held for $t\geq0$ and $x\in S.$

Property 1) is easily justified.  Applying to (14) the theorem on
the jump of the co-normal derivative of a single-layer potential
[6] or, more precisely, its version for the integrals over
$[t,\infty)$ instead of the ones over $[0,t]$, we get the
relations
$$
\frac{\partial V_\lambda(t,x\pm)}{\partial
N(x)}=\mp\psi(t,x)\pm\lambda V_\lambda(t,x)+
\int^\infty_td\tau\int_S \frac{\partial\widetilde
V(\tau-t,x,y)}{\partial N(x)}(\psi(\tau,y)-\lambda
V_\lambda(\tau,y)) d\sigma_y
$$ valid for $t>0$ and $x\in S$.

Taking into account (9),(10), we arrive at formula (15).

Obviously, the restriction of the function $V_\lambda(t,x)$ on
$[0,\infty)\times S$ is a solution to the equation (8) in the
required class. We now show that there is no more than one such a
solution. Assume that $f_1(t,x)$ and $f_2(t,x)$ are two solutions
from this class. Put $\hat f(t,x)=f_1(t,x)-f_2(t,x)$. Then there
exists a function $H\hat f$, and the relation
$$
\lambda \hat f(t,x)= \left[ \frac{1+q(x)}{2}\frac{\partial H\hat
f(t,x+)}{\partial N(x)}- \frac{1-q(x)}{2}\frac{\partial H\hat
f(t,x-)}{\partial N(x)} \right]\tag16
$$
fulfilled for $t>0,\ x \in S$. Choose $T_0>0$ such that for  all
$t\geq T_0,\ x\in S,\ f_1(t,x)=f_2(t,x)=0$. If $\underset {t\in
[0,T_0],\ x\in S}\to{\inf}H\hat f(t,x)=\beta<0$, then according to
the maximum principle there exists a point $(t_0,x_0)\in
[0,T_0)\times S$ such that $H \hat f(t_0,x_0)=\hat
f(t_0,x_0)=\beta$. This implies the inequalities $$ \frac{\partial
H\hat f(t_0,x_0-)}{\partial N(x_0)}\leq 0, \frac{\partial H\hat
f(t_0,x_0+)}{\partial N(x_0)}\geq 0.\tag17
$$

From (17) we have that the right-hand side of (16) is non-negative
at the point $(t_0,x_0)$. But this contradicts the assertion  that
$\hat f(t_0,x_0)<0$. Thus $\underset {t\in [0,T_0],\ x\in
S}\to{\inf} \hat f(t,x)=$ \break $\underset {t\in [0,T_0],\ x\in
\Bbb R^d}\to{\inf}H \hat f(t,x)\geq0.$ We can get that $\underset
{t\in [0,T_0],\ x\in S}\to{\sup} \hat f(t,x)\leq0$ in the same
manner. So, $\hat f(t,x)\equiv0$ on $(t,x)\in [0,\infty)\times S$.
This completes the proof for $r$ being equal to $0$. In the case
of non-negative $r$ we get the assertion from the previous one
arguing as in [4], pp. 194-196.
\enddemo

\proclaim {Proposition 4}
Let $r$ and $q$ be given continuous
real-valued functions on $S$ such that $r$ is bounded and
non-negative, $|q|\leq 1.$ Then the uniqueness theorem is valid
for the boundary process.
\endproclaim

\demo
 { Proof} We follow the
proof of Theorem 5.2 in [4].  The martingale property can be
easily extended to continuous functions on $[0,\infty)\times S$
having compact supports  with respect to $t$ in the manner of
Remark 2. Given $x\in S,$ for any measure ${\Bbb R}_x$ on $\Cal
D([0,+\infty),[0,+\infty)\times S)$ being a solution to the
submartingale problem starting from $x$ we have the relation
$$
\Bbb E^{{\Bbb R}_x}[f(\tau(\theta), y(\theta))]= f(0,x)+\Bbb
E^{{\Bbb R}_x} \left[
\int^\theta_0(\widetilde{K}f)(\tau(u),y(u)))du \right].
$$
Performing the Laplace transformation, we get,
for $\lambda>0,$ the equality
$$
\int^\infty_0e^{-\lambda u}\Bbb E^{{\Bbb R}_x} [\lambda
f(\tau(u),x(\tau(u)))- (\widetilde{K}f)(\tau(u),
x(\tau(u)))]du=f(0,x).
$$
Then, for $\psi\in C^\infty_0([0,\infty)\times S)$, the integral
$\int^\infty_0e^{-\lambda u}\Bbb E^{{\Bbb R}_x}\psi(\tau(u),
x(\tau(u))) du$ is uniquely determined, provided the equation
$\lambda f-\widetilde{K}f=\psi$ has the unique solution for each
\break  $\psi\in C^\infty_0([0,\infty)\times S)$. But this
condition is true because of the Lemma. The assertion of the
Proposition follows in the way of Corollary 6.2.4 in [9].
\enddemo
\head {3. The main result} \endhead
 \proclaim { Theorem} Let $S$ be a
closed bounded surface in ${\Bbb R}^d$ which belongs to the class
$H^{2+\delta}$ for some $\delta\in (0,1)$, $q$ and $r$ be given
continuous functions on $S$ taking values in $[-1, 1]$ and
$[0,+\infty)$, respectively, and $r$ is bounded. For $y\in \Bbb
{R}^d$, let $b(y)$ be a symmetrical $d\times d$-matrix satisfying
the conditions $J$. Then, for each $x\in{\Bbb R}^d$, there exists
a unique solution to the submartingale problem.
\endproclaim

\demo
 { Proof}  This assertion follows from
Proposition 5 by the arguments of Theorem 4.2 in [4].
\enddemo


\Refs \ref \no1 \by O.V.Aryasova, M.I. Portenko \paper One example
of a random change of time that transforms a generalized
 diffusion process into an ordinary one
\jour Theory Stochast. Process \vol 13(29) \pages  \issue 3 \yr
2007
\endref
\ref \no2 \by N.I. Portenko \book Generalized diffusion processes
\publ Naukova Dumka \publaddr Kiev \yr1982 \transl English transl.
\publ Amer. Math. Soc., Providence, RI \yr 1990
\endref
\ref \no3
 \by  E.B. Dynkin
\book Markov processes \publaddr Moscow \publ Fizmatgiz \yr 1963.
\transl English transl \vol I, II \publ New York: Acad.Press.
Berlin: Springer \yr 1965
\endref
\ref\no4 \by D.W.Strook, S.R.S.Varadhan \paper Diffusion processes
with boundary conditions \jour  Comm. Pure Appl. Math. \vol 23
\pages147--225 \yr1971
\endref
\ref \no5 \by O.V.Aryasova, M.I. Portenko \paper One class of
multidimensional stochastic differential equations having no
property of weak uniqueness of a solution \jour Theory Stochast.
Process \vol 11(27) \pages 14--28 \yr 2005 \issue 3--4
\endref
\ref \no6 \by O. Ladyzhenskaya, V. Solonnikov, and N. Ural'tseva
\book Linear and Quasilinear Equations of Parabolic Type \publaddr
Moscow \publ Nauka \yr 1967 \transl Translations of Mathematical
Monographs, American Mathematical Society  \vol 23 \yr 1968
\publaddr Providence, RI
\endref
\ref \no7 \by S.V.Anulova \paper Diffusion processes with singular
characteristics \jour Int. Symp. Stochast. Dif\-fe\-rent. \break
Equat.: Abstract. Vilnius \yr1978 \pages7--11
\endref
\ref \no8 \by A. Friedman \paper Partial differential equations of
parabolic type \jour Prentice--Hall, Inc., Englewood Cliffs, N.J.
\yr1964
\endref
\ref \no9 \by D.W.Strook, S.R.S.Varadhan \book Multidimensional
diffusion processes \publaddr Berlin \publ Springer \yr 1979
\endref
\endRefs
\enddocument